# Points rationnels dans leur fibre : compléments à un théorème de Poonen

Laurent Moret-Bailly[*][†][‡]

**Abstract**

Let $f : X \to S$ be a morphism of integral schemes, which is surjective and of finite presentation. When $S$ is of finite type over a field or over $\mathbb{Z}$, of positive dimension, Poonen has shown that there is a point $x \in X$ whose residue field $\kappa(x)$ is purely inseparable over $\kappa(f(x))$. In this paper, we extend this result in several ways. First we prove that we can take $x$ such that $f(x)$ is a codimension 1 point of $S$; if $S$ is smooth over a fiel $k$, we can require $\kappa(f(x))$ to be separable over $k$. In another direction, we prove that similar results hold for other schemes $S$, such as noetherian integral schemes of dimension $\geq 2$.

**Résumé**

Soit $f : X \to S$ un morphisme surjectif et de présentation finie de schémas intègres. Lorsque $S$ est de type fini sur un corps ou sur $\mathbb{Z}$, et de dimension $> 0$, Poonen a montré qu'il existe un point $x \in X$ dont le corps résiduel $\kappa(x)$ est une extension radicielle de $\kappa(f(x))$. Dans ce travail, on prolonge ce résultat de plusieurs façons. D'abord, on peut choisir $x$ tel que $f(x)$ soit un point de codimension 1 de $S$ ; si $S$ est lisse sur un corps $k$, on peut exiger que $\kappa(f(x))$ soit séparable sur $k$. Dans une autre direction, on montre des résultats analogues pour d'autres classes de schémas $S$, par exemple les schémas noethériens intègres de dimension $\geq 2$.

*Classification AMS :* 14 A15, 14 D99.

## 1 Introduction

La source de ce travail est le théorème suivant, dû à Poonen :

**Théorème 1.1.** [7, Theorem 1] *Soit $k$ un corps et soit $f : X \to S$ un morphisme dominant de $k$-variétés, avec $\dim(S) \geq 1$. Alors il existe un point $x$ de $X$ dont le corps résiduel $\kappa(x)$ est une extension radicielle de $\kappa(f(x))$.* □

En particulier, si $k$ est parfait, il existe un point $s$ de $S$ dont la fibre $f^{-1}(s)$ admet un point $\kappa(s)$-rationnel. Il en est de même d'ailleurs si $f$ est lisse (*loc. cit.*, remarque suivant la preuve du théorème). Il y a un énoncé analogue pour les morphismes de $\mathbb{Z}$-schémas de type fini [7, Theorem 4].

Le but de cet article est double. D'une part, sous les hypothèses du théorème de Poonen, on apporte les compléments suivants :

[*]Univ Rennes, CNRS, IRMAR - UMR 6625, F-35000 Rennes, France
[†]laurent.moret-bailly[AT]univ-rennes1.fr

- on peut trouver un point $x$, comme dans l'énoncé, tel que $f(x)$ soit un point de codimension 1 de $S$ ;
- si $S$ est lisse sur $k$, on peut exiger en outre que $\kappa(f(x))$ soit une extension séparable de $k$ ;
- si $X$ est lisse sur $k$, il existe un point fermé $x$ de $X$ tel que $\kappa(x) = \kappa(f(x))$ et que $\kappa(x)$ soit séparable sur $k$. (Ce résultat figure aussi dans [6] : voir le lemme 14 de *loc. cit.* et les quelques lignes qui le précèdent).

D'autre part, on cherche à généraliser ce type de résultats à d'autres schémas de base $S$. À cet effet, on introduit la notion suivante. Soit $S$ un schéma, que pour simplifier nous supposerons affine et intègre dans l'introduction. On dira que $S$ est un *schéma de Poonen* (voir la définition 3.1 pour les schémas généraux) si pour tout $S$-schéma $f : X \to S$ dominant et de présentation finie, il existe un point $x$ de $X$ tel que $\kappa(x)$ soit radiciel sur $\kappa(f(x))$. (Ainsi, le théorème 1.1 affirme que tout schéma intègre de type fini sur un corps, de dimension $\geq 1$, est un schéma de Poonen.) On montre alors (théorème 3.2) que si $f : X \to S$ est comme ci-dessus, et si $S$ est un schéma de Poonen, il existe un diagramme commutatif

$$\begin{array}{ccc} X' & \xhookrightarrow{j} & X \\ \downarrow{f'} & & \downarrow{f} \\ S' & \xhookrightarrow{i} & S \end{array} \qquad \text{(D)}$$

où $S'$ est non vide, $i$ et $j$ sont des immersions, $f'$ est fini, surjectif et radiciel, et $S'$ est *localement défini par une équation* dans $S$ : ainsi, le phénomène de « codimension 1 » rencontré plus haut est général.

On cherche ensuite des conditions suffisantes pour qu'un schéma (toujours supposé affine et intègre, disons $S = \mathrm{Spec}\,(A)$) soit un schéma de Poonen. Il s'avère qu'il suffit de considérer des morphismes $f : X \to S$ de la forme $X = \mathrm{Spec}\,(A[Y]/(F))$ où $F$ est un polynôme d'une variable non constant, à coefficients dans $A$. Il y a plusieurs énoncés de ce type, le plus simple étant le suivant (condition (P2) de 4.1.5) : $S$ est un schéma de Poonen si, pour tout $F \in A[Y]$ non constant, il existe $b \in A$ tel que $F(b) \notin A^\times$. (Cette condition apparaît dans [7], sous une forme un peu différente). À partir de ce critère et de ses variantes, on retrouve le théorème de Poonen, mais aussi d'autres exemples : ainsi, tout schéma noethérien intègre de dimension $\geq 2$ est un schéma de Poonen, ainsi que tout schéma intègre de dimension $\geq 2$ dont les anneaux locaux sont factoriels.

## 1.2 Plan de l'article

Dans un souci de cohérence, on a pris le parti de donner des démonstrations indépendantes de [7], même lorsque le recours à [7] aurait permis d'abréger le texte.

Le §2 contient des définitions et des énoncés de base. Au §3, on définit les schémas de Poonen et l'on prouve l'existence du diagramme (D) considéré plus haut (théorème 3.2). Au §4, on donne des conditions suffisantes pour qu'un schéma soit un schéma de Poonen, en partant de critères polynomiaux tels que la condition (P2), et l'on en déduit les principaux exemples connus : schémas de type fini sur un corps ou sur $\mathbb{Z}$ (4.2.3), « courbes » de Fargues-Fontaine, certains schémas intègres de dimension $\geq 2$ (4.3). On donne aussi quelques exemples de schémas de Poonen stricts (4.4). Les compléments au théorème de Poonen évoqués au début sont démontrés au §5 ; on y établit en fait une version relative de ces résultats (théorème 5.1), où le schéma $S$ est supposé lisse sur un schéma de base $B$ quelconque. Il est à noter que le §5 n'utilise que de façon marginale les résultats des sections précédentes ; avec quelques aménagements, on obtient ainsi une démonstration du théorème 5.1 indépendante de ces sections.

## 1.3 Remerciements

Je remercie le rapporteur pour ses remarques et pour ses questions (qui m'ont notamment conduit aux résultats des sections 4 et 5), et Jean-Louis Colliot-Thélène qui m'a signalé l'article [6].

## 1.4 Notations et conventions

Les anneaux sont commutatifs et unitaires. Si $A$ est un anneau et $t$ un élément de $A$, le localisé $(t^{\mathbb{N}})^{-1}A$ sera noté $A_t$ ; on notera $D(t)$ l'ouvert $\mathrm{Spec}\,(A_t)$ de $\mathrm{Spec}\,(A)$ et $V(t)$ le fermé $\mathrm{Spec}\,(A/tA)$.

Si $s$ est un point d'un schéma $S$, l'anneau local de $S$ en $s$ et le corps résiduel de $s$ seront notés respectivement $\mathscr{O}_{S,s}$ et $\kappa(s)$.

*Conventions pour les flèches des diagrammes :*

| | | | | | |
|---|---|---|---|---|---|
| $\hookrightarrow$ | = | monomorphisme | $\hookrightarrow\!\!\circ\!\!\to$ | = | immersion ouverte |
| $\twoheadrightarrow$ | = | morphisme surjectif | $\hookrightarrow\!\!/\!\!\to$ | = | immersion fermée |
| $\rightarrowtail$ | = | morphisme radiciel | $\hookrightarrow\!\!\phi\!\!\to$ | = | immersion |

# 2 Quasi-relèvements

**Définition 2.1.** *Soit $f : X \to S$ un morphisme de schémas, et soit $T \xrightarrow{i} S$ un $S$-schéma. Un* relèvement *de $T$ à $X$ (ou à $f$) est un $S$-morphisme $j : T \to X$ tel que $f \circ j = i$.*

*Un* quasi-relèvement *de $T$ à $X$ est un diagramme commutatif*

$$\begin{array}{ccc} \widetilde{T} & \xrightarrow{j} & X \\ \downarrow{\rho} & & \downarrow{f} \\ T & \xrightarrow{i} & S \end{array}$$

*où $\rho$ est* fini, surjectif et radiciel.

*Nous dirons que $T$ est (quasi-)relevable à $X$ s'il admet au moins un (quasi-)relèvement.*

Dans ce travail, nous n'utiliserons ce langage que lorsque $i$ est un monomorphisme, et plus précisément dans deux cas :
– $T = \mathrm{Spec}\,(\kappa)$ est un point de $S$, $i$ est l'inclusion, et $\widetilde{T}$ est un point de $X$ au-dessus de $T$, à corps résiduel fini radiciel sur $\kappa$ ;
– $i$ est une immersion.

**Définition 2.2.** *Si $f : X \to S$ est un morphisme de schémas, nous poserons*

$$\begin{aligned}
\widetilde{\mathrm{QRel}}(f) &:= \{x \in X \mid \text{l'extension } \kappa(f(x)) \subset \kappa(x) \text{ est radicielle}\} \\
\widetilde{\mathrm{Rel}}(f) &:= \{x \in X \mid \text{l'extension } \kappa(f(x)) \subset \kappa(x) \text{ est triviale}\} \\
\mathrm{QRel}(f) &:= f(\widetilde{\mathrm{QRel}}(f)) \subset S \\
\mathrm{Rel}(f) &:= f(\widetilde{\mathrm{Rel}}(f)) \subset S.
\end{aligned}$$

Ainsi, $\mathrm{Rel}(f)$ (resp. $\mathrm{QRel}(f)$) est l'ensemble des points de $S$ qui sont relevables (resp. quasi-relevables) à $X$.

**Définition 2.2.1.** *Une* immersion principale *est un morphisme $u : X \to Y$ de schémas qui se décompose sous la forme $X \overset{i}{\hookrightarrow} U \overset{j}{\hookrightarrow} Y$ où $j$ est une immersion ouverte et $i$ une immersion fermée localement définie par une équation*[1].

*Si $X$ est un sous-schéma de $Y$ et $u$ l'inclusion, on dira aussi que $X$ est un* sous-schéma principal *de $Y$, et que son espace sous-jacent est un* sous-espace principal *de $|Y|$.*

**Proposition 2.2.2.** *Soit $f : X \to S$ un morphisme étale de schémas.*

*(1) On a $\widetilde{\mathrm{QRel}}(f) = \widetilde{\mathrm{Rel}}(f)$ et $\mathrm{QRel}(f) = \mathrm{Rel}(f)$.*

*(2) Soit $x \in \widetilde{\mathrm{Rel}}(f)$ et soit $s = f(x)$. Il existe un sous-schéma* principal *$S' \overset{i}{\hookrightarrow} S$ de $S$ contenant $s$, et un relèvement $j$ de $S'$ à $X$ qui est une immersion principale et vérifie $j(s) = x$.*

*Preuve.* L'assertion (1) est immédiate (et il suffit que $f$ soit non ramifié) puisque toutes les extensions résiduelles de $f$ sont séparables. Montrons (2). On peut remplacer $X$ et $S$ par des voisinages convenables de $x$ et $s$ ; d'après le théorème de structure locale des algèbres étales [8, V, théorème 1], nous pouvons donc supposer que $S = \mathrm{Spec}(A)$ et $X = \mathrm{Spec}(B)$ où $B = (A[Y]/(F))_h$, $Y$ étant une indéterminée et $h$ un élément de $A[Y]/(F)$. Puisque $\kappa(s) \overset{\sim}{\to} \kappa(x)$, l'image de $Y$ dans $\kappa(x)$ définit un élément $\alpha$ de $\kappa(s)$ (vérifiant $F(\alpha) = 0$ et $h(\alpha) \neq 0$) ; quitte à localiser $A$, on peut relever $\alpha$ en $a \in A$. Soit alors $S' := V(F(a))$ le sous-schéma fermé principal défini par $F(a)$ : on obtient $j : S' \to X$ à partir du morphisme d'algèbres $B \to A/(F(a))$ envoyant la classe $y$ de $Y$ sur la classe de $a$ : ce morphisme identifie $A/(F(a))$ à $B/(y-a)$, ce qui prouve que $j$ est une immersion principale. Les inclusions $i(S') \subset \mathrm{Rel}(f)$ et $j(S') \subset \widetilde{\mathrm{Rel}}(f)$ sont triviales. $\square$

**Remarque 2.2.3.** En termes géométriques, la construction de la démonstration se décrit ainsi : quitte à localiser, on a une $S$-immersion principale $u$ de $X$ dans la droite affine $\mathbb{A}_S^1$ (de coordonnée $Y$) ; d'autre part le point $x$ est image d'une section de $\mathbb{A}_S^1$ au-dessus de $\mathrm{Spec}(\kappa(s))$, que l'on prolonge en une section $\sigma$ au-dessus de $S$ (donnée par $Y = a$). Le schéma $S'$ est le produit fibré $X \times_{u, \mathbb{A}_S^1, \sigma} S$ ; le fait que $i$ et $j$ soient des immersions principales résulte par changement de base de la même propriété pour $u$ et $\sigma$.

**Définition 2.2.4.** *Un morphisme de schémas $f : X \to S$ sera dit* quasi-étale *s'il se décompose en $X \overset{\rho}{\to} X_0 \overset{\pi}{\to} S$ où $\pi$ est étale et où $\rho$ est entier, surjectif et radiciel (ce qui équivaut à dire que $\rho$ est un homéomorphisme universel, cf. [5, 18.12.11]).*

**Proposition 2.2.5.** *Soit $f : X \to S$ un morphisme quasi-étale, décomposé en $X \overset{\rho}{\to} X_0 \overset{\pi}{\to} S$ comme dans 2.2.4.*

---
[1]On n'impose aucune condition à l'équation : elle peut être nulle, par exemple.

(1) *On a* $\mathrm{QRel}(f) = \mathrm{QRel}(\pi) = \mathrm{Rel}(\pi)$.

(2) *Soit* $x \in \widetilde{\mathrm{QRel}}(f)$ *et soit* $s = f(x)$. *Il existe un sous-schéma principal* $S' \overset{i}{\hookrightarrow} S$ *contenant* $s$, *et un quasi-relèvement* $S' \overset{f'}{\twoheadleftarrow} \widetilde{S}' \overset{j}{\hookrightarrow} X$ *de* $S'$ *à* $X$, *où* $\widetilde{S}'$ *contient* $x$ *et où* $j$ *est une immersion principale.*

*Preuve.* On laisse (1) au lecteur (la seconde égalité est un rappel de 2.2.2). Pour montrer (2), on applique 2.2.2 (2) au morphisme $\pi$, ce qui donne des immersions principales $i : S' \hookrightarrow S$ et $j_0 : S' \hookrightarrow X_0$ avec $\pi \circ j_0 = i$. Il suffit alors de poser $X' := X \times_{\rho, X_0, j_0} S'$.  □

**Remarque 2.2.6.** On voit ainsi que si $f : X \to S$ est quasi-étale, alors $\mathrm{QRel}(f)$ (resp. $\widetilde{\mathrm{QRel}}(f)$) est réunion de sous-espaces principaux de $|S|$ (resp. de $|X|$), et est en particulier *ind-constructible* dans $S$ (resp. dans $X$) ; pour cette notion, voir [4, §7.2].

On peut vérifier que cette dernière propriété est valable pour tout morphisme $f : X \to S$ localement de présentation finie.

La proposition suivante est utile pour ramener « génériquement » l'étude de $\mathrm{QRel}(f)$ au cas quasi-étale. Noter que l'hypothèse que $S$ est réduit est en général inoffensive puisque $\mathrm{QRel}(f) = \mathrm{QRel}(f_{\mathrm{red}})$.

**Proposition 2.3.** *Soient* $S$ *un schéma réduit*, $f : X \to S$ *un morphisme localement de présentation finie*, $\eta$ *un point* maximal *de* $S$, $\xi$ *un point* fermé *de* $X_\eta = f^{-1}(\eta)$.

*Il existe un diagramme commutatif comme ci-contre, dans lequel* $S_0$ *est un voisinage ouvert affine de* $\eta$, $v : X_0 \to X$ *est une immersion, et* $f_0$ *est quasi-étale, fini, localement libre et surjectif.*

$$\begin{array}{ccc} \xi \in X_0 & \overset{v}{\hookrightarrow} & X \\ \downarrow f_0 & & \downarrow f \\ \eta \in S_0 & \overset{u}{\hookrightarrow} & S \end{array}$$

*De plus, si l'extension* $\kappa(\eta) \subset \kappa(\xi)$ *est séparable, on peut prendre* $f_0$ *étale.*

*Preuve.* Puisque $X_\eta$ est un $\kappa(\eta)$-schéma de type fini, $\kappa(\xi)$ est une extension finie de $\kappa(\eta)$, et admet donc une sous-extension séparable maximale $L$ ; l'extension $\kappa(\xi)/L$ est finie radicielle.

Puisque $S$ est réduit et $\eta$ maximal, $\kappa(\eta)$ coïncide avec $\mathscr{O}_{S,\eta}$, et donc $\mathrm{Spec}(\kappa(\eta))$ est limite projective de la famille filtrante des voisinages affines de $\eta$ ; par suite il y a un tel voisinage $S_0 = \mathrm{Spec}(A_0)$ et un morphisme $B_0 \to C_0$ de $A_0$-algèbres de présentation finie, induisant $L \to \kappa(\xi)$ par le changement de base $A_0 \to \kappa(\eta)$. Quitte à localiser, $B_0$ est finie étale sur $A_0$ et $C_0$ est libre et radicielle sur $B_0$.

Comme $f$ est localement de présentation finie, les résultats de [5, §8] montrent que, quitte à localiser à nouveau, l'immersion fermée naturelle $\mathrm{Spec}(\kappa(\xi)) \to X_\eta$ se prolonge en une immersion $\mathrm{Spec}(C_0) \to X$ avec les propriétés voulues. La dernière assertion de l'énoncé vient du fait que si $\kappa(\eta) \subset \kappa(\xi)$ est séparable, on a $L = \kappa(\xi)$ et l'on peut prendre $C_0 = B_0$.  □

## 3 Schémas de Poonen

**Définition 3.1.** *Un schéma* $S$ *est un* schéma de Poonen *(resp. un* schéma de Poonen strict*) s'il vérifie la condition suivante :*

*Pour tout ouvert non vide* $U \subset S$ *et tout* $U$-*schéma* $f : X \to U$ *de présentation finie et dominant, on a* $\mathrm{QRel}(f) \neq \emptyset$ *(resp.* $\mathrm{Rel}(f) \neq \emptyset$*).*

*Un anneau $A$ est un* anneau de Poonen (strict) *si* $\operatorname{Spec}(A)$ *est un schéma de Poonen (strict).*

Dans ce travail, on s'intéressera surtout aux schémas et anneaux de Poonen ; voir le §4.4 pour quelques exemples de schémas de Poonen stricts.

Dans le cas d'un schéma noethérien, on peut simplifier la définition 3.1 en oubliant l'ouvert $U$. Plus généralement :

**Proposition 3.1.1.** *Soit $S$ un schéma.*

*(1) Si $S$ est un schéma de Poonen (resp. un schéma de Poonen strict), alors pour tout $g : Y \to S$ de présentation finie et dominant, $\operatorname{QRel}(g)$ (resp. $\operatorname{Rel}(g)$) est dense dans $S$.*

*(2) La réciproque est vraie si $S$ est quasi-séparé et n'a qu'un nombre fini de composantes irréductibles.*

*Preuve.* L'assertion (1) est immédiate. Sous les hypothèses de (2), si $U$ est un ouvert affine non vide de $S$ et $f : X \to U$ un morphisme dominant de présentation finie, il existe, quitte à restreindre $U$, un ouvert affine $V$ disjoint de $U$ et tel que $U \cup V$ soit dense dans $S$ (considérer les points maximaux de $S$ n'appartenant pas à $U$). Il suffit alors d'appliquer l'hypothèse au morphisme $g : X \amalg V \to S$, qui est bien de présentation finie si $X$ est quasi-séparé, les immersions de $U$ et $V$ dans $X$ étant quasi-compactes. $\square$

### 3.1.2 Exemples et remarques

**3.1.2.1.** J'ignore si l'assertion 3.1.1 (2) est vraie sans hypothèse sur $S$.

**3.1.2.2.** Si un morphisme $f : X \to U$ est de présentation finie et dominant, son image contient un ouvert dense de $U$ ; il en résulte que l'on peut, dans la définition, se limiter aux morphismes $f$ qui sont surjectifs, et exiger que $\operatorname{QRel}(f)$ (resp. $\operatorname{Rel}(f)$) soit dense dans $U$. On peut aussi se limiter aux morphismes $f : X \to U$ où $X$ et $U$ sont affines.

**3.1.2.3.** $S$ est un schéma de Poonen (strict) si et seulement si $S_{\mathrm{red}}$ en est un.

**3.1.2.4.** La somme $\coprod_{i \in I} S_i$ d'une famille $(S_i)_{i \in I}$ de schémas est un schéma de Poonen (strict) si et seulement si chacun des $S_i$ en est un.

**3.1.2.5.** Un corps est un anneau de Poonen (strict) si et seulement si il est séparablement clos (algébriquement clos).

**3.1.2.6.** Soit $S'$ un ouvert de $S$ : si $S$ est un schéma de Poonen (strict), il en est de même de $S'$. La réciproque est vraie si $S'$ est dense dans $S$ ; plus généralement :

**3.1.2.7.** Soit $u : S' \to S$ un morphisme *dominant*. On suppose que pour tout $s' \in S'$, l'extension résiduelle $\kappa(u(s')) \subset \kappa(s')$ est radicielle (resp. triviale). Alors, si $S'$ est un schéma de Poonen (resp. un schéma de Poonen strict), il en est de même de $S$.
En effet, traitons le cas des schémas de Poonen (le cas strict est analogue). Soit $f : X \to U \subset S$ comme dans la définition 3.1 : alors $U \times_S S' \to U$ est encore dominant à extensions résiduelles radicielles, et $U \times_S S'$ contient un ouvert affine non vide $V$ de $S'$. Le morphisme induit $X \times_S V \to V$ est dominant ; on a donc $\operatorname{QRel}(X \times_S V / V) \neq \emptyset$, et le fait que les extensions résiduelles de $V \to U$ soient radicielles entraîne que $u(\operatorname{QRel}(X \times_S V / V)) \subset \operatorname{QRel}(X/U)$.

**3.1.2.8.** Soit $u : S' \to S$ un morphisme de schémas. Dans chacun des cas suivants, $u$ vérifie les hypothèses de 3.1.2.7., donc $S$ est un schéma de Poonen (resp. un schéma de Poonen strict) si $S'$ en est un :

(i) $u$ est radiciel et dominant (resp. un monomorphisme dominant) ;
(ii) $S$ est irréductible et $S' = \mathrm{Spec}\,(\mathscr{O}_{S,s})$, où $s$ est un point de $S$ ;
(iii) $S'$ est la somme des composantes irréductibles de $S$.

**3.1.2.9.** Soit $f : S \to T$ un morphisme de schémas. Supposons que pour tout $t \in T$ la fibre $S_t := f^{-1}(t)$ soit un schéma de Poonen (strict). Alors $S$ est un schéma de Poonen (strict) : il suffit en effet d'appliquer 3.1.2.7. avec $S' = \coprod_{t \in T} S_t$.

**3.1.2.10.** Soit $S$ un schéma de Poonen (strict), et soit $u : S' \to S$ un $S$-schéma de présentation finie. On suppose que pour tout ouvert non vide $V$ de $S'$, $u(V)$ contient un ouvert non vide de $S$ (de façon équivalente, $u$ envoie tout point maximal de $S'$ sur un point maximal de $S$). Alors $S'$ est un schéma de Poonen (strict).

**Théorème 3.2.** *Soit $S$ un schéma de Poonen non vide, et soit $f : X \to S$ un morphisme dominant de présentation finie.*

*(1) Il existe un sous-schéma* principal *non vide $S'$ de $S$ qui est quasi-relevable à $X$, et même relevable à $X$ si $f$ est lisse.*

*(2) On suppose $S$ localement noethérien. Alors $\mathrm{QRel}\,(f)$ contient un point de codimension $\leq 1$ de $S$ ; il en est de même de $\mathrm{Rel}\,(f)$ si $f$ est lisse.*

*(3) On suppose que $S$ est un schéma de Jacobson [4, (6.4.1)]. Alors $\mathrm{QRel}\,(f)$ contient un point fermé de $S$ ; il en est de même de $\mathrm{Rel}\,(f)$ si $f$ est lisse.*

*Preuve.* On peut supposer $S$ réduit. Les assertions (2) et (3) sont des conséquences immédiates de (1). Montrons (1). L'image de $f$ contient un point maximal $\eta$ de $S$ ; choisissons un point fermé $\xi$ de $X_\eta$, avec $\kappa(\xi)$ séparable sur $\kappa(\eta)$ si $f$ est lisse. Appliquant 2.3, on voit que quitte à remplacer $S$ par un voisinage de $\eta$, on peut supposer $f$ quasi-étale, et même étale dans le cas lisse. Comme $S$ est un schéma de Poonen, il existe un point $s \in S$ qui est dans $\mathrm{QRel}\,(f)$, et donc dans $\mathrm{Rel}\,(f)$ dans le cas étale. Il suffit alors d'appliquer 2.2.5 (2), et 2.2.2 (2) dans le cas étale. □

**Remarque 3.2.1.** Si l'on suppose dans 3.2 que $S$ est un schéma de Poonen strict, on ne peut pas pour autant conclure (en dehors du cas lisse) que $\mathrm{Rel}\,(f)$ contient des points de codimension 1. Par exemple, soit $k$ un corps parfait de caractéristique $p > 0$, et considérons $S = \mathbb{A}_k^2$. Nous verrons plus loin (4.4.3) que $S$ est un schéma de Poonen strict. Soit $f : S \to S$ le morphisme de Frobenius relatif à $k$ : on voit facilement que $\mathrm{Rel}\,(f)$ est l'ensemble des points fermés de $S$.

# 4 Exemples de schémas de Poonen

## 4.1 Cas des schémas affines intègres : critères polynomiaux

Le point de départ de nos exemples est l'énoncé suivant, conséquence facile de 2.3 :

**Proposition 4.1.1.** *Soit $S$ un schéma. Les conditions suivantes sont équivalentes :*

(i) $S$ est un schéma de Poonen ;

(ii) pour tout ouvert affine non vide $U$ de $S$ et tout morphisme fini, surjectif et quasi-étale $f : X \to U$, on a $\operatorname{QRel}(f) \neq \emptyset$ ;

(iii) pour tout ouvert affine non vide $U$ de $S$ et tout morphisme fini étale surjectif $f : X \to U$, on a $\operatorname{QRel}(f) \neq \emptyset$ (et donc $\operatorname{Rel}(f) \neq \emptyset$ d'après 2.2.2 (1)).

*Preuve.* Les implications (i)⇒(ii)⇒(iii) sont triviales, et 2.2.5 (1) montre que (iii) implique (ii). Supposons (ii), et montrons (i). Soit donc $f : X \to U$ dominant de présentation finie, où $U$ est un ouvert affine non vide de $S$. Choisissons un point maximal $\eta$ de $U$ contenu dans l'image de $f$, et un point fermé $\xi$ de $f^{-1}(\eta)$. Appliquant 2.3, on trouve, quitte à restreindre $U$, un $U$-schéma fini, quasi-étale et surjectif $f' : X' \to U$ et une $U$-immersion fermée $X' \hookrightarrow X$. Il est clair que $\operatorname{QRel}(f') \subset \operatorname{QRel}(f)$ et l'hypothèse (ii) assure que $\operatorname{QRel}(f') \neq \emptyset$, donc $\operatorname{QRel}(f) \neq \emptyset$. □

**Notation 4.1.2.** Soit $A$ un anneau intègre, de spectre $S$, de corps des fractions $K$, et soit $Y$ une indéterminée. Pour tout $F \in A[Y]$, nous poserons

$$X_{A,F} := \operatorname{Spec}(A[Y]/(F)).$$

Si aucune confusion n'est à craindre, nous écrirons $X_F$ plutôt que $X_{A,F}$. Nous n'utiliserons ces notations que lorsque $F$ n'est pas constant ; dans ce cas, $X_F$ est un $S$-schéma dominant et génériquement fini. Nous dirons que $F$ est *séparable* si son image dans $K[Y]$ est un polynôme séparable, c'est-à-dire si $X_F$ est génériquement étale sur $S$.

Dans ce qui suit, nous nous intéresserons aux ensembles de la forme $\operatorname{Rel}(X_F/S)$, et secondairement $\operatorname{QRel}(X_F/S)$, le but étant de déterminer si $A$ est un anneau de Poonen. On laisse au lecteur la démonstration du lemme suivant :

**Lemme 4.1.3.** *Soit $A$ un anneau intègre, de corps des fractions $K$ et de spectre $S$, et soit $F \in A[Y]$ non constant.*

(1) *Soit $\mathfrak{p}$ un idéal premier de $A$. Les conditions suivantes sont équivalentes :*

   (i) $\mathfrak{p} \in \operatorname{Rel}(X_F/S)$ ;
   (ii) *l'image de $F$ dans $\kappa(\mathfrak{p})[Y]$ a une racine dans $\kappa(\mathfrak{p})$ ;*
   (iii) *il existe $g \in A \smallsetminus \mathfrak{p}$ et $b \in A_g$ tels que $F(b) \in \mathfrak{p} A_{\mathfrak{p}}$.*

(2) *Les conditions suivantes sont équivalentes :*

   (i) $\operatorname{Rel}(X_F/S)$ *est dense dans $S$ ;*
   (ii) *pour tout $g \in A \smallsetminus \{0\}$, il existe $u \in A$ et $b \in A_{ug}$ tels que $F(b) \notin A_{ug}^{\times}$.*
   (iii) *pour tout $g \in A \smallsetminus \{0\}$, il existe $u \in A$, $b \in A_{ug}$ et $\mathfrak{p} \in D(ug)$ tels que $F(b) \in \mathfrak{p} A_{\mathfrak{p}}$.* □

La proposition 4.1.1 conduit au critère suivant :

**Proposition 4.1.4.** *Soit $A$ un anneau intègre, de spectre $S$. Pour que $A$ soit un anneau de Poonen, il faut et il suffit que pour tout $F \in A[Y]$ séparable non constant, $\operatorname{Rel}(X_F/S)$ soit dense dans $S$.*

*Preuve.* La nécessité est triviale. Montrons la suffisance en utilisant le critère 4.1.1 (iii). Soient donc $U \subset S$ un ouvert affine non vide, $f : X \twoheadrightarrow U$ un morphisme fini étale surjectif, et montrons que $\mathrm{Rel}\,(X/U) \neq \emptyset$. La fibre générique $X_K$ de $f$ (où $K := \mathrm{Frac}(A)$) est le spectre d'un produit fini d'extensions finies séparables de $K$ ; quitte à restreindre encore $U$ et à remplacer $X$ par un ouvert fermé, on peut supposer $X$ intègre. De plus, d'après le théorème de l'élément primitif (et quitte à restreindre $U$ à nouveau), on peut supposer que $U = D(g)$ pour $g \neq 0$ dans $A$, et que $X = X_{A_g,F}$ où $F \in A_g[Y]$ est non constant et séparable sur $K$. On peut écrire $F = g^{-N} F_1$ où $N \in \mathbb{N}$ et $F_1 \in A[Y]$ ; il est clair que $X$ est $U$-isomorphe à $X_{A_g,F_1}$. Appliquant l'hypothèse à $F_1$, on conclut que $\mathrm{Rel}\,(X_{A,F_1}/S) \cap U \neq \emptyset$ et donc que $\mathrm{Rel}\,(X/S) \cap U \neq \emptyset$. □

**Notation 4.1.5.** Pour un anneau intègre $A$, de spectre $S$, nous considérerons les conditions

(P1)    Pour tout $F \in A[Y]$ non constant, $\mathrm{Rel}\,(X_F/S)$ est dense dans $S$.

(P2)    Pour tout $F \in A[Y]$ non constant, il existe $b \in A$ tel que $F(b) \notin A^\times$.

**Remarques 4.1.6.** Gardons les notations de 4.1.5.

(1) Si $A$ est un corps, (P1) (resp. (P2)) signifie que $A$ est algébriquement clos.

(2) (P2) implique que le radical de $A$ est nul : considérer $F(Y) = 1 + aY$, avec $a \neq 0$.

(3) Dans la démonstration du lemme 3 de [7], Poonen mentionne la condition suivante : pour tout $F \in A[Y]$ non constant et tout $g \neq 0$ dans $A$, il existe $b \in A$ tel que $F(b) \notin A_g^\times$. Cette condition est en fait équivalente à (P2), comme le montre la preuve de l'implication (c) de la proposition 4.1.8 plus bas.

Pour étudier les liens entre ces conditions, le lemme élémentaire qui suit nous sera utile ; la preuve en est laissée au lecteur.

**Lemme 4.1.7.** *Soient $A$ un anneau, $g$ et $b$ deux éléments de $A$ tels que les images de $b$ dans $A/(g)$ et dans $A_g$ soient inversibles. Alors $b \in A^\times$.* □

**Proposition 4.1.8.** *Soient $A$ un anneau intègre et $g \in A$ non nul.*

(1) *On a les implications :*

$$
\begin{array}{ccccc}
A \text{ vérifie (P2)} & \stackrel{(a)}{\Longrightarrow} & A \text{ vérifie (P1)} & \stackrel{(b)}{\Longrightarrow} & A \text{ est un anneau de Poonen} \\
\Downarrow (c) & & \Updownarrow (d) & & \\
A_g \text{ vérifie (P2)} & \stackrel{(a)}{\Longrightarrow} & A_g \text{ vérifie (P1)} & &
\end{array}
$$

(2) *Si $A$ est un anneau de Jacobson, (a) et (c) sont des équivalences.*

*Preuve.* (1) L'implication (b) résulte de 4.1.4 (la différence étant que, dans (P1), le polynôme $F$ considéré n'est pas supposé séparable).

L'équivalence (d) est facile et laissée au lecteur. Montrons (c). Supposons donc (P2) pour $A$, soit $F \in A_g[Y]$ non constant, et cherchons $b \in A_g$ tel que $F(b) \notin A_g^\times$. Pour $N$ entier convenable, $F_1 := g^N F \in A[Y]$ et, pour tout $b \in A_g$, on a l'équivalence $F(b) \in A_g^\times \Leftrightarrow F_1(b) \in A_g^\times$. On peut donc remplacer $F$ par $F_1$ et supposer $F = \sum_{i=0}^d a_i Y^i \in A[Y]$, avec $d > 0$ et $a_d \neq 0$. Si $a_0 = 0$, alors $b = 0$ convient. Sinon, considérons le polynôme

$$H(Y) := a_0^{-1} F(a_0 g Y) = 1 + a_1 g Y + \cdots + a_d a_0^{d-1} g^d Y^d.$$

Comme $H$ n'est pas constant (car $a_d a_0^{d-1} g^d \neq 0$), l'hypothèse (P2) fournit $b' \in A$ tel que $H(b') \notin A^\times$. Comme d'autre part $H(b') \in 1 + gA$, on déduit du lemme 4.1.7 que $H(b') \notin A_g^\times$, et donc que $F(a_0 g b') = a_0 H(b') \notin A_g^\times$, d'où la conclusion avec $b = a_0 g b'$.

Montrons (a). Supposons que $A$ vérifie (P2), et soit $F \in A[Y]$ non constant. Pour voir que $\mathrm{Rel}\,(X_F/S)$ est dense dans $S$, appliquons le critère 4.1.3 (2). Soit donc $0 \neq g \in A$ : on cherche $u \in A$ et $b \in A_{ug}$ tels que $F(b) \notin A_{ug}^\times$. Or nous venons de voir que $A_g$ vérifie (P2), donc il existe même un tel couple $(u, b)$ avec $u = 1$.

(2) Supposons que $A$ soit un anneau de Jacobson. Vu les implications déjà établies, il suffit de voir que (P1) implique (P2) pour $A$. Soit $F \in A[Y]$ non constant. (P1) implique l'existence de $u \in A$ non nul et de $b \in A_u$ tel que $F(b) \notin A_u^\times$. Soit $\mathfrak{n}$ un idéal maximal de $A_u$ contenant $F(b)$ : puisque $A$ est un anneau de Jacobson, $\mathfrak{m} := \mathfrak{n} \cap A$ est un idéal maximal de $A$ de même corps résiduel $\kappa$. Il existe $b' \in A$ dont l'image dans $\kappa$ est celle de $b$ : on a alors $F(b') \in \mathfrak{m}$, donc $F(b') \notin A^\times$. $\square$

**Remarques 4.1.9.** (1) En général, les implications (a), (b) et (c) sont strictes. C'est immédiat pour (b) : un corps séparablement clos, mais non parfait, est un anneau de Poonen qui ne vérifie pas (P1).

Pour (a), prenons $A$ *local* noethérien intègre de dimension $\geq 2$. La remarque 4.1.6 (2) montre que $A$ ne vérifie pas (P2), mais nous verrons plus bas (4.3.1) que $A$ vérifie (P1). Le même exemple vaut pour (c), car pour tout $g \neq 0$ dans l'idéal maximal de $A$, l'anneau $A_g$ vérifie (P2) (remarque 4.3.2).

(2) On a introduit la condition (P2) parce qu'elle est, lorsqu'elle s'applique, facile à vérifier puisqu'elle n'utilise pas de localisation : on va le voir dans la section 4.2.

(3) Le statut logique de (P2) semble aussi plus simple. Si l'on fixe le degré de $F$ dans (P2), la condition s'exprime par une formule du premier ordre du langage des anneaux. Par suite la classe des anneaux intègres vérifiant (P2) est élémentaire ; en particulier, si $A$ vérifie (P2), tout « modèle non standard » de $A$ vérifie (P2). J'ignore si la classe des anneaux de Poonen (ou celle des anneaux vérifiant (P1)) est élémentaire.

## 4.2 Exemples valués et conséquences

Dans cette section, nous appellerons *valeur absolue* sur un corps $K$ une application $|.| : K \to \overline{\Gamma}$, où :
– $(\Gamma, \times, \leq)$ est un groupe abélien totalement ordonné, et $\overline{\Gamma}$ est le monoïde $\Gamma \cup \{0\}$, ordonné par $0 < \Gamma$ ;
– ou bien $\overline{\Gamma} = \mathbb{R}_{\geq 0}$ et $|.|$ est une valeur absolue (au sens habituel), ou bien $|.|$ est une valuation sur $K$, notée multiplicativement (avec notamment $|0| = 0$).

**Théorème 4.2.1.** *Soit $A$ un anneau intègre, dont le corps des fractions $K$ est muni d'une valeur absolue $|.| : K \to \overline{\Gamma}$. On suppose que le sous-monoïde $|A|$ de $\overline{\Gamma}$ n'est pas majoré, mais que $|A^\times|$ l'est.*

*Alors $A$ vérifie la condition* (P2) *de 4.1.5 et est donc un anneau de Poonen.*

Preuve. Soit $F(Y) = a_d Y^d + G(Y) \in A[Y]$, avec $d > 0$, $a_d \neq 0$, et $\deg(G) < d$. Les hypothèses entraînent que si $b \in A$ avec $|b|$ suffisamment grand, on a :
– dans le cas ultramétrique, $|a_d b^d| > |G(b)|$ et donc $|F(b)| = |a_d| |b|^d$ ;
– pour une valeur absolue classique, $|a_d b^d| > 2 |G(b)|$ et donc $|F(b)| \geq \frac{1}{2} |a_d| |b|^d$.
On en déduit que $|F(A)|$ n'est pas majoré et donc que $F(A) \not\subset A^\times$. $\square$

**Corollaire 4.2.2.**

(1) [7, Lemma 3] $\mathbb{Z}$ *est un anneau de Poonen.*

(2) [7, Theorem 1] *Si $k$ est un corps, l'anneau de polynômes $k[t]$ est un anneau de Poonen.*

*Preuve.* Il suffit d'appliquer 4.2.1 avec, pour $\mathbb{Z}$, la valeur absolue ordinaire, et pour $k[t]$ la valeur absolue de Gauss $P \mapsto 2^{-\deg(P)}$. □

**Théorème 4.2.3.**

(1) [7, Theorem 1] *Soit $k$ un corps. Tout $k$-schéma de type fini sans point isolé est un schéma de Poonen.*

(2) [7, Theorem 4] *Tout $\mathbb{Z}$-schéma de type fini sans point isolé est un schéma de Poonen.*

*Preuve.* Dans les deux cas, on se ramène (en utilisant les remarques 3.1.2) au cas d'un schéma $S$ affine intègre de dimension $\geq 1$. Dans le cas (1), il existe un $k$-morphisme dominant $S \to \mathbb{A}^1_k$, de sorte que d'après 3.1.2.10., l'assertion résulte du cas de $\mathbb{A}^1_k = \operatorname{Spec}(k[t])$. Dans le cas (2), ou bien $S$ domine $\operatorname{Spec}(\mathbb{Z})$ et l'on conclut par 4.2.2 (1), ou bien $S$ est de caractéristique $p > 0$, c'est un $\mathbb{F}_p$-schéma de type fini et on applique (1). □

Voici une généralisation du cas des anneaux de polynômes, qui s'applique notamment aux « courbes » construites dans [3][2] :

**Théorème 4.2.4.** *Soit $\overline{S}$ un schéma intègre, de corps des fonctions $K$, muni d'un point fermé $\infty$ vérifiant les propriétés suivantes :*

(i) *l'anneau local $\mathscr{O} := \mathscr{O}_{\overline{S}, \infty}$ est un anneau de valuation ;*

(ii) $S := \overline{S} \smallsetminus \{\infty\}$ *est un schéma affine non vide, d'anneau noté $A$ ;*

(iii) $\Gamma(\overline{S}, \mathscr{O}_{\overline{S}}) = A \cap \mathscr{O}$ *est un corps.*

*Alors $A$ vérifie* (P2) *et est donc un anneau de Poonen (et par suite, d'après 3.1.2.6., $\overline{S}$ est un schéma de Poonen).*

*Preuve.* Notons $\mathfrak{m}$ l'idéal maximal de $\mathscr{O}$ et $|.| : K \to \overline{\Gamma} = K/\mathscr{O}^\times$ la valuation associée, qui est non triviale puisque $S \neq \emptyset$. L'hypothèse (iii) entraîne que $A \cap \mathfrak{m} = \{0\}$ ; pour tout $a \in A$ non nul, on a donc $|a| \geq 1$, et par suite $|a| = 1$ si $a \in A^\times$. De plus $|A|$ n'est pas majoré dans $\Gamma$ : en effet, pour tout $z = a/b \in K$ (avec $a$ et $b \neq 0$ dans $A$) on a $|a| = |b||z| \geq |z|$ donc $|A|$ est cofinal dans $\Gamma$, qui n'est pas majoré puisque $|.|$ est non triviale. Les conditions de 4.2.1 sont donc satisfaites, ce qui établit le résultat. □

**Remarque 4.2.5.** Soit $A$ un anneau principal qui n'est pas un corps. Pour que $A$ soit un anneau de Poonen, il faut qu'il ait une infinité d'idéaux maximaux (sinon le point générique est un ouvert de $\operatorname{Spec}(A)$ qui n'est pas un schéma de Poonen). Cette condition n'est pas suffisante. Par exemple, soient $P$ l'ensemble des nombre premiers, $\Sigma = \{p \in P \mid p = 2 \text{ ou } p \equiv 1 \pmod 4\}$, et $\Sigma' = \{p \in P \mid p \equiv -1 \pmod 4\}$. Alors $\Sigma$ et $\Sigma'$ sont de même densité $1/2$, $A_1 := \Sigma^{-1}\mathbb{Z}$ n'est pas un anneau de Poonen (considérer $\operatorname{Spec}(A_1[Y]/(Y^2 + 1))$), alors que $A_2 := \Sigma'^{-1}\mathbb{Z}$ vérifie (P2) : pour $F \in A_2[Y]$ non constant, soit $\alpha$ une racine de $F$ dans $\mathbb{C}$ et soit $K := \mathbb{Q}(\alpha, i)$ : si $p$ est un nombre premier décomposé dans $K$, on a $p \notin \Sigma'$ (donc $(p) \in \operatorname{Spec}(A_2)$) et $(p)$ est relevable à $A_2[Y]/(F)$.

---

[2]Cependant beaucoup de ces courbes (celles du chapitre 6 de [3]) sont *trivialement* des schémas de Poonen puisque les corps résiduels de leurs points fermés sont algébriquement clos.

## 4.3  Schémas de dimension ≥ 2

**Théorème 4.3.1.** *Soit $A$ un anneau local intègre de dimension $\geq 2$, d'idéal maximal $\mathfrak{m}$. On suppose vérifiée la condition suivante :*

(*) *Pour tout $x \in \mathfrak{m}$ non nul, $\mathrm{Spec}\,(A/xA)$ n'a qu'un nombre fini de composantes irréductibles, et elles sont de codimension $1$ dans $\mathrm{Spec}\,(A)$.*

*Alors, $A$ vérifie la condition* (P1) *de 4.1.5 et est donc un anneau de Poonen.*

*Preuve.* Soit $F = \sum_{i=0}^{d} a_i Y^i \in A[Y]$, avec $d > 0$ et $a_d \neq 0$, et montrons que $\mathrm{Rel}\,(X_F/S)$ est dense dans $S := \mathrm{Spec}\,(A)$. Pour cela, fixons $g \in A$ non nul (que l'on peut supposer dans $\mathfrak{m}$). Nous allons exhiber $b \in A_{ga_d}$ et $\mathfrak{q} \in D(ga_d)$ tels que $F(b) \in \mathfrak{q}A_{\mathfrak{q}}$, ce qui établira le théorème d'après 4.1.3 (2).

Soient $\mathfrak{p}_i$ ($1 \leq i \leq n$) les idéaux premiers de $A$ correspondant aux composantes de $V(ga_d)$. Le lemme d'évitement des idéaux premiers [2, II, §1, n° 1, prop. 2] montre qu'il existe $c \in \mathfrak{m} \smallsetminus \cup_{i=1}^{n}\mathfrak{p}_i$. Prenons $b := c/(ga_d^2) \in A_{ga_d}$. On trouve

$$F(b) = a_d^{1-2d} g^{-d} \left[ c^d + \sum_{i=0}^{d-1} a_i a_d^{2(d-i)-1} g^{d-i} c^i \right] = a_d^{1-2d} g^{-d} \left( c^d + g a_d \alpha \right)$$

où $\alpha \in A$. L'élément $\gamma := c^d + g a_d \alpha$ appartient à $\mathfrak{m}$ (comme $c$ et $g$) et n'est dans aucun des $\mathfrak{p}_i$, vu le choix de $c$. Si $\mathfrak{q} \in \mathrm{Spec}\,(A)$ est un point maximal de $V(\gamma)$, alors $\mathfrak{q}$ est de hauteur $1$ d'après (*), et donc $ga_d \notin \mathfrak{q}$ (sinon $\mathfrak{q}$ serait l'un des $\mathfrak{p}_i$, ce qui est exclu car il contient $\gamma$). Autrement dit, $\mathfrak{q} \in D(ga_d)$ et $\gamma \in \mathfrak{q}A_{\mathfrak{q}}$, et donc $F(b) \in \mathfrak{q}A_{\mathfrak{q}}$, le facteur $a_d^{1-2d} g^{-d}$ étant dans $A_{ga_d}$. □

**Remarque 4.3.2.** Sous les hypothèses de 4.3.1, soit $g \in \mathfrak{m} \smallsetminus \{0\}$. Supposons que $A_g$ soit un anneau de Jacobson : alors, puisqu'il vérifie (P1), il vérifie aussi (P2) d'après 4.1.8 (2). Cela se produit notamment si $A$ est *noethérien* [5, (10.5.8)], et aussi lorsque $\dim(A) = 2$ : dans ce dernier cas, $A_g$ est un anneau intègre de dimension $1$ qui (d'après le lemme d'évitement) admet une infinité d'idéaux maximaux.

**Corollaire 4.3.3.** *Soit $S$ un schéma localement noethérien, dont toutes les composantes irréductibles sont de dimension $\geq 2$. Alors $S$ est un schéma de Poonen.*

*Preuve.* Utilisant 3.1.2.8. et 3.1.2.3., nous sommes ramenés au cas où $S$ est le spectre d'un anneau local noethérien intègre de dimension $\geq 2$, qui est un cas particulier de 4.3.1. □

**Corollaire 4.3.4.** *Soit $A$ un anneau de Krull [2, VII, §1, n° 3, définition 3] de dimension $\geq 2$. Alors $A$ est un anneau de Poonen.*

*En particulier, tout anneau factoriel de dimension $\geq 2$ est un anneau de Poonen.*

*Preuve.* Soit $\mathfrak{p}$ un idéal premier de $A$ de hauteur $\geq 2$. Alors $A_{\mathfrak{p}}$ est un anneau de Krull [2, VII, §1, n° 3, prop. 6] auquel le théorème 4.3.1 s'applique ; c'est donc un anneau de Poonen. Par suite $A$ est un anneau de Poonen d'après 3.1.2.8.. □

## 4.4  Exemples de schémas de Poonen stricts

**Proposition 4.4.1.** *Soit $S$ un schéma de Poonen intègre, dont le corps de fonctions $K$ est parfait. Alors $S$ est un schéma de Poonen strict.*

*Preuve.* Soit $f : X \to S$ de présentation finie et dominant. Pour voir que $\mathrm{Rel}\,(f) \neq \emptyset$, on peut supposer $X$ réduit. Puisque $K$ est parfait, l'ouvert de lissité de $f$ est non vide, ce qui nous ramène au cas où $f$ est lisse, qui résulte de 3.2 (1). □

**Proposition 4.4.2.** *Soit $S$ un schéma de Poonen. On suppose que $S$ est un schéma de Jacobson et que, pour tout point fermé $s$ de $S$, le corps résiduel $\kappa(s)$ est parfait. Alors $S$ est un schéma de Poonen strict.*

*Preuve.* C'est une conséquence immédiate de 3.2 (3). □

**Corollaire 4.4.3.**

(1) *Soit $k$ un corps* parfait. *Tout $k$-schéma de type fini sans point isolé est un schéma de Poonen strict.*

(2) *Tout $\mathbb{Z}$-schéma de type fini sans point isolé est un schéma de Poonen strict.* □

# 5 Schémas lisses sur une base

**Théorème 5.1.** *Soit $X \xrightarrow{f} S \xrightarrow{g} B$ un diagramme de schémas non vides. On suppose que :*

(i) *$g$ est lisse, surjectif, à fibres de dimension $\geq 1$ ;*

(ii) *$f$ est de présentation finie et dominant.*

(1) *Il existe un sous-schéma non vide $S'$ de $S$, lisse de codimension relative 1 sur $B$, quasi-relevable à $X$, et relevable à $X$ si $f$ est lisse.*

(2) *Si $X$ est lisse sur $B$, il existe un sous-schéma non vide $S'$ de $S$, étale sur $B$, et relevable à $X$.*

## 5.2 Preuve de 5.1 : réduction au cas d'un corps de base

Les arguments étant classiques, on laisse les détails au lecteur. En remplaçant $X$, $S$ et $B$ par des ouverts convenables, on peut supposer qu'ils sont affines et que la dimension des fibres de $g$ est constante. Posant $B = \mathrm{Spec}\,(k)$, on remarque que puisque $f$ et $g$ sont de présentation finie, il existe un sous-anneau $k_0$ de $k$ de type fini sur $\mathbb{Z}$ et un diagramme $X_0 \xrightarrow{f_0} S_0 \xrightarrow{g_0} B_0 := \mathrm{Spec}\,(k_0)$ induisant $X \xrightarrow{f} S \xrightarrow{g} B$ par changement de base et ayant les mêmes propriétés : on peut donc supposer $B$ noethérien. Si $\eta$ est un point maximal de $B$ et $B_1 := \mathrm{Spec}\,(\mathscr{O}_{B,\eta})$, il suffit de démontrer le théorème pour le diagramme $X_1 \xrightarrow{f_1} S_1 \xrightarrow{g_1} B_1$ obtenu par le changement de base $B_1 \to B$ : un sous-schéma $S_1' \subset S_1$ de présentation finie est induit par un sous-schéma de $S$, qui aura les propriétés voulues.

Ceci nous ramène au cas où $B = \mathrm{Spec}\,(k)$ est local artinien, de corps résiduel $\kappa$. Par changement de base à $B_2 := \mathrm{Spec}\,(\kappa)$, on obtient un diagramme $X_2 \xrightarrow{f_2} S_2 \xrightarrow{g_2} B_2$. Si $S_2' \subset S_2$ est un sous-schéma lisse sur $B_2$, il se relève (au moins localement) en un sous-schéma $S' \subset S$ lisse sur $B$, de même dimension relative ; si $S_2'$ est quasi-relevable à $X_2$, alors $S'$ est trivialement quasi-relevable à $X$ (tout schéma fini, radiciel et surjectif sur $S_2'$ est aussi fini, radiciel et surjectif sur $S'$) ; si $f$ est lisse, tout relèvement de $S_2'$ à $X_2$ se prolonge en un relèvement de $S'$ à $X$, au moins si $S'$ est affine.

En résumé, il suffit d'établir le théorème lorsque $B$ est le spectre d'un corps.

## 5.3 Preuve de 5.1 : le cas d'un corps de base

On suppose désormais que $B = \mathrm{Spec}\,(k)$, où $k$ est un corps ; $S$ est un $k$-schéma affine lisse, purement de dimension $d > 0$, et $f : X \to S$ est de type fini et dominant.

Supposons d'abord $k$ parfait. Alors le théorème 5.1 est une conséquence immédiate de 3.2 et du fait que, d'après 4.2.3, $S$ est un schéma de Poonen : il suffit de remarquer, pour (1), que tout sous-schéma réduit de $S$ est génériquement lisse sur $k$ et, pour (2), que tout point fermé de $S$ est étale sur $k$.

La suite de la démonstration s'applique dès que $k$ est *infini*, ce que nous supposons désormais ; il est à noter que les arguments qui suivent ne font pas appel à 3.2 ni à 4.2.3, qui sont donc redémontrés dans ce cas. (Pour le cas d'un corps fini, voir la remarque 5.3.1.1.)

### 5.3.1 Le cas où $\dim(S) = 1$

Comme toujours, on suppose $X$ et $S$ affines, et $S$ irréductible de dimension 1. On peut remplacer $X$ par une courbe irréductible de $X$, dominant $S$, lisse si $X$ l'est, et étale sur $S$ si $f$ est lisse. Ainsi nous pouvons supposer, en préservant les hypothèses dans chaque cas, que $X$ et $S$ sont des courbes affines.

Montrons d'abord (2). On suppose donc $X$ lisse sur $k$, et l'on cherche un point fermé de $S$, séparable sur $k$ et relevable à $X$. Fixons un $k$-morphisme dominant $\pi : S \to \mathbb{A}^1_k$. Il suffit de démontrer le résultat avec $f$ remplacé par $\pi \circ f$. On supposera donc que $S = \mathbb{A}^1_k = \mathrm{Spec}\,(k[t])$.

Le corps de fonctions $K := k(X)$ est une extension séparable de $k$, de degré de transcendance 1 ; en particulier le degré d'imperfection de $K$ sur $k$ est égal à 1 (en d'autres termes, $[K : kK^p] = p$, où $p$ est l'exposant caractéristique de $k$) ; *a fortiori*, le degré d'imperfection de $K$ sur $k(S)$ est $\leq 1$. Par suite (voir [1, V, §13, exercice 3] ou [9, §40]) l'extension finie $k(X)/k(S)$ est *engendrée par un élément*. Ceci nous ramène au cas où $X \subset \mathbb{A}^2_k$ est un ouvert de $X_{k[t],F}$ (notation de 4.1.2) où $F \in A[Y] = k[t,Y]$ est non constant en $Y$. Il existe une droite $D \subset \mathbb{A}^2_k$ d'équation $Y = \lambda t + \mu$ ($\lambda, \mu \in k$) qui coupe $X$ transversalement en au moins un point $x$ : en effet, les droites ayant cette propriété forment un ouvert $U$ du plan des $(\lambda, \mu)$ ; on a $U \neq \emptyset$ car $X$ est lisse, donc $U(k) \neq \emptyset$ puisque $k$ est infini. Puisque $x \in D$, on a $x \in \widetilde{\mathrm{Rel}}\,(X/S)$, et puisque l'intersection est transverse, $\kappa(x)$ est séparable sur $k$, ce qui prouve (2).

Montrons (1). Le cas où $f$ est étale résulte de (2) ; en général, on peut supposer que $f$ se décompose en $X \xrightarrow{g} X_0 \xrightarrow{f_0} S$ où $f_0$ est étale et $g$ radiciel. Il suffit d'appliquer le cas étale à $f_0$ et de remarquer que $\mathrm{QRel}\,(X/S) = \mathrm{Rel}\,(X_0/S)$.

**Remarque 5.3.1.1.** L'hypothèse que $k$ est infini a servi pour l'existence de la droite transverse $D$ dans la preuve de (2). Si $k$ est fini, on peut procéder comme suit : on peut supposer que $X$ est l'ouvert $g \neq 0$ de $X_{k[t],F}$, où $g \in k[t]$. Si $F(t,0) = 0$, alors $X$ contient l'ouvert $g \neq 0$ de l'axe $Y = 0$, et la conclusion est triviale. Sinon, un calcul élémentaire montre que, pour $N$ assez grand, $F(t, g^{-N})$ n'est pas inversible dans $k[t]_g$ ; un idéal maximal qui le contient fournit un point fermé de $\mathrm{Rel}\,(X/S)$, automatiquement séparable puisque $k$ est parfait. Cet argument élimine le recours à 3.2 et à 4.2.3.

### 5.3.2 Le cas général

Commençons par montrer (2). Si $X$ est lisse, il existe une courbe lisse irréductible $C \subset X$ telle que $f(C)$ soit une courbe dans $S$, automatiquement lisse (au moins génériquement). Il suffit d'appliquer le cas de dimension 1 au morphisme induit $C \to f(C)$.

Montrons enfin (1). Ceci revient à trouver un point $s \in \mathrm{QRel}(f)$, de codimension 1 dans $S$, séparable sur $k$, et appartenant à $\mathrm{Rel}(f)$ si $f$ est lisse.

On peut supposer que $\dim(S) = d > 1$, et qu'il existe un morphisme lisse $\pi : S \to \mathbb{A}_k^{d-1}$. Soit $K$ le corps des fonctions rationnelles de $\mathbb{A}_k^{d-1}$. Le changement de base $\mathrm{Spec}(K) \to \mathbb{A}_k^{d-1}$ donne un $K$-morphisme $f_K : X_K \to S_K$ vérifiant encore les hypothèses du théorème (avec $f_K$ lisse si $f$ l'est, et $S_K$ lisse sur $K$) et avec $\dim(S_K) = 1$. Le cas de dimension 1 fournit un point fermé $s_K \in \mathrm{QRel}(f_K)$, séparable sur $K$ et dans $\mathrm{Rel}(f_K)$ si $f$ est lisse. Ce point (vu comme point de $S$) répond à la question : il est clair que $\mathrm{QRel}(f_K) \subset \mathrm{QRel}(f)$ et $\mathrm{Rel}(f_K) \subset \mathrm{Rel}(f)$ puisque $X_K \to X$ et $S_K \to S$ sont des monomorphismes ; d'autre part $s$ est de codimension 1 dans $S$ (on a $\mathscr{O}_{S,s} = \mathscr{O}_{S_K,s}$), et enfin l'extension $k \subset \kappa(s)$ est composée des extensions séparables $k \subset K \subset \kappa(s)$, donc est séparable. $\square$

# Références